\documentclass{amsart}%
\usepackage{amsmath}
\usepackage{graphicx}
\usepackage{amscd}%
\usepackage{amsfonts}%
\usepackage{amssymb}
\providecommand{\U}[1]{\protect\rule{.1in}{.1in}}
\newtheorem{theorem}{Theorem}
\theoremstyle{plain}

\newtheorem{axiom}{Axiom}

\newtheorem{definition}{Definition}

\newtheorem{remark}{Remark}

\numberwithin{equation}{section}

\begin{document}
\title[\textsf{Introduction to: Non-Associaitive Finite Invertible Loops}]{\textsf{Introduction to: Non-Associaitive Finite Invertible Loops}}
\author{\textsf{Raoul E. Cawagas}}
\address{Raoul E. Cawagas, SciTech R\&D Center, Polytechnic University of the
Philippines, Sta. Mesa, Manila}
\email{raoulec@yahoo.com}
\thanks{This paper is part of a series of papers on NAFILs (\emph{Non-Associaitive
Finite Invertible Loops}) to appear in forthcoming issues of this journal (PUPJST).}

\begin{abstract}
\emph{Non-associative finite invertible loops (NAFIL)} are loops whose every
element has a unique two-sided inverse. Not much is known about the class of
NAFIL loops which includes the familiar IP (Inverse Property), Moufang, and
Bol loops. Our studies have shown that they are involved in such diverse
fields as combinatorics, finite geometry, quasigroups and related systems,
Cayley algebras, as well as in theoretical physics. This paper presents an
introduction to the class of NAFIL loops as the starting point for the
development of the theory of these interesting structures.

\end{abstract}
\maketitle

\section{\textbf{Introduction}}

Many associative algebraic structures like \emph{groups} have been studied
extensively by numerous mathematicians [2]. It appears, however, that not much
is known in the current literature about non-associative group-like structures
other than special classes of \emph{quasigroups }and \emph{loops} because the
theory of these structures [13]. ``... is a fairly young discipline which
takes its roots from geometry, algebra and combinatorics.''

In 1981, we became interested in studying a class of non-associative finite
loops in which every element has a unique inverse [3]. Several structures of
this type have been intensively studied like IP (Inverse Property), Moufang,
and Bol loops. However, there are many interesting loops in this class other
than these. A literature search has shown that there is no standard term for
this class of loops. Because of this, we have decided to introduce the acronym
\emph{\textbf{NAFIL}} which stands for: \emph{\textbf{non-associative finite
invertible loop}}.

This paper presents an introduction to the class of \emph{NAFIL }loops in an
attempt to start the development of the theory of these structures. It is
based on the results of various studies started in 1981 that led the way to a
two-year research project supported by the National Research Council of the
Philippines (NRCP)\footnote{Supported by the National Research Council of the
Philippines under NRCP Project B-88 and Project B-95: \emph{Development of the
Theory of Finite Pseudogroups.}}. This project [6] which was conducted at the
SciTech R\&D Center of the Polytechnic University of the Philippines from 1996
to 1998 is still being actively pursued at present.

To our knowledge, NAFIL loops have not been systematically studied as a
specific class of finite loops. Moreover, very little is known about their
applications in other fields of mathematics and in theoretical physics. For
instance, our studies have shown that these loops are involved in the theories
of non-associative real algebras and rings. Thus, among other things, we have
shown [4] that the basis vectors of all non-associative $2^{r}$-dimensional
\emph{Cayley algebras \ }form NAFIL loops of order $2^{r+1}$. These loops
define real algebras that constitute an infinite family which includes all
known alternative real division algebras as members. One of these algebras
(the \emph{Cayley numbers }or \emph{octonions}) has also been applied by
physicists to the quantum theory of elementary particles in connection with
quark models and in string theory [12].

In this paper, we shall introduce the class of NAFILs by defining the concept
of the \textbf{\emph{non-associative finite invertible loop (NAFIL)}} and
presenting some of its fundamental properties. As a sequel to this
introductory paper, we shall deal with the structure of NAFILs, develop tools
and methods (algorithms) for their construction and analysis, identify
fundamental problems to be solved, and present some of the important results
we have so far obtained from our researches in an attempt to lay the
foundations of the\textbf{\emph{\ Theory of Non-Associative Finite Invertible
Loops}}\emph{. }

\section{\textbf{The Finite Invertible Loop}}

To fully appreciate what the invertible loop is, it is helpful to recall the
idea of an abstract mathematical system. Such a system essentially consists
of: a non-empty set $S$ of distinct elements, at least one binary operation
$\star$, an equivalence relation $=$, a set of axioms, as well as a set of
definitions and theorems and is usually denoted by $\left(  S,\star\right)  $.
The heart of the system is the set of axioms from which all the theorems are
derived. Most algebraic systems like groups, rings, and fields satisfy some or
all of the following axioms or postulates:

\begin{axiom}
For all $a,b\in S,$ $a\star b\in S$. \emph{(Closure axiom)}
\end{axiom}

\begin{axiom}
There exists a unique element $e\in S,$ called the \emph{identity}, such that
$e\star a=a\star e=a$ for all $a\in S$. \emph{(Identity axiom)}
\end{axiom}

\begin{axiom}
Given an identity element $e\in S,$ for every $a\in S$ there exists a unique
element $a^{-1}\in S,$ called its \emph{inverse}, such that $a\star
a^{-1}=a^{-1}\star a=e$. \emph{(Inverse axiom)}
\end{axiom}

\begin{axiom}
For every $a,b\in$ $S$ there exists unique $x$, $y\in$ $S$ such that $a\star
x=b$ and $y\star a=b$. \emph{(Unique Solution axiom)}
\end{axiom}

\begin{axiom}
For every $a,b\in S$, $a\star b=b\star a$\emph{. (Commutative axiom)}
\end{axiom}

\begin{axiom}
For all $a,b,c\in S$, $(a\star b)\star c=a\star(b\star c)$. \emph{(Associative
axiom)}
\end{axiom}

The simplest algebraic system is the \emph{groupoid}; it is only required to
satisfy A1. This is a trivial system and not much can be said about it. A
groupoid that also satisfies A4 is called a \emph{quasigroup}; and a
quasigroup that satisfies A2 is called a \emph{loop} [9, 13].

One of the most important systems is the \emph{group} which can be defined as
any system satisfying axioms A1, A2, A3, A4, and A6. Hence, the group is a
loop that also satisfies A3 and A6. Alternatively, we can say that the group
is a quasigroup that satisfies A2, A3, and A6. It must be noted that the
axioms A1, A2, A3, A4, A5, and A6 are independent. (Of course A3 is meaningful
only if A2 is assumed.) Any algebraic system that satisfies A5 (Commutative
axiom) is called \emph{abelian} (or \emph{commutative}).

A loop that also satisfies A3 is called an \emph{invertible loop}. Such a loop
can be either associative or non-associative.

\begin{remark}
The class of \textbf{loops} forms a \textbf{variety}\footnote{A \emph{variety
}is\emph{\ }a class of algebraic structures of the same signature satisfying a
given set of identities which is closed under the taking of homomorphic images
(H), subalgebras (S), and direct products (P). Thus, the class of invertible
loops also forms a variety. The class of NAFIL loops, however, does not form a
variety because it is not closed under H and P; a NAFIL can have homomorphic
images and subloops that are groups (which are not NAFILs).} [9, 10]. defined
by the set of operations $\Omega=\left\{  e,^{-1},\star\right\}  $ (where $e,$
$^{-1},$ and $\star$ are nullary, unary, and binary operations, respectively),
a set of identities $I=\{e\star x=x,$ $x\star e=x\}$ (where $e$ is a unique
identitiy element), and the equations: $ax=b,$ $ya=b$.
\end{remark}

\begin{definition}
An\emph{\ }\textbf{invertible loop} $\left(  \mathcal{L},\star\right)  $ is a
non-empty set $\mathcal{L}$ with a binary operation $\star$ satisfying axioms
A1, A2, A3, and A4. If it also satisfies A5, it is called \textbf{abelian} and
if $\mathcal{L}$ has a finite number of elements, it is called \textbf{finite}.
\end{definition}

It is clear from Definition 1 that \emph{an invertible loop is a loop in which
every element has a unique (or two-sided) inverse.} Hence the class of all
invertible loops belongs to the \emph{variety} of loops and is defined by the
set of operations $\Omega=\left\{  e,^{-1},\star\right\}  ,$ the equations
$ax=b,$ $ya=b,$ and the identities:
\[
e\star x=x,\text{ }x\star e=x\quad x\star x^{-1}=e,\text{ }x^{-1}\star x=e
\]
where $e$ is the identity element and $x^{-1}$ is the unique inverse of $x.$
Thus groups are also invertible loops. The term \emph{invertible loop}
therefore applies to both associative and non-associative types. Moreover, it
can be shown that the class of invertible loops also forms a variety of loops.

\section{\textbf{The Non-Associative Finite Invertible Loop (NAFIL)}}

Having defined what the invertible loop is, let us now focus our attention on
the class of \textbf{\emph{non-associative finite invertible loops }}whose
abstract theory we are attempting to develop.

\begin{definition}
A finite invertible loop that is non-associative is called a \textbf{NAFIL
(non-associative finite invertible loop) }while one that is associative is
called a \textbf{group}.
\end{definition}

To emphasize the fact that the acronym \textbf{\emph{NAFIL}} refers to a loop,
we shall allow ourselves some grammatical freedom by often calling it a
\textbf{\emph{NAFIL loop}}. Although NAFIL loops are invertible loops, the
class of NAFILs does not form a variety.

The class of NAFIL loops includes, among others, the following:

\begin{itemize}
\item Loops with inverse properties: IP (Inverse Property), LIP/RIP
(Left/Right Inverse Property), CIP (Crossed Inverse Property), etc.

\item Moufang loops and Bol loops

\item Plain loops (anti-associative)
\end{itemize}

\begin{remark}
\emph{We note that a distinction is sometimes made between the terms
}non-associative\textbf{\emph{\ }}\emph{\ and }not-associative\emph{.\ A
system is ``non-associative'' if A6 (associative axiom) is not assumed to hold
while it is ``not-associative'' if A6 is required to be not satisfied by the
system as a whole. The NAFIL, as defined above, is therefore not-associative
in this sense. However, unless otherwise indicated, the term non-associative
will henceforth be taken to mean not-associative when referring to NAFIL
loops.}
\end{remark}

In this paper, our main concern will be on finite invertible loops that are
\emph{not-associative }(NAFIL loops). This does not imply, however, that A6
does not hold at all in such a loop; it may be satisfied in a limited way
within the system. Since all finite algebraic systems we will consider satisfy
A1 (Closure axiom), they are essentially \emph{groupoids}. Henceforth, we
shall therefore use \textbf{\emph{groupoid}} as a generic term for finite
quasigroups, loops, groups, and related structures. For such finite systems,
it is both convenient and useful to define them in \emph{constructive} terms.

A finite algebraic system or groupoid like the invertible loop $\left(
\mathcal{L},\star\right)  $ is completely defined by its multiplication or
\emph{Cayley table}. Such a table is a listing of the $n^{2}$ possible binary
products, $\ell_{i}\star\ell_{j},$ of its $n$ elements and it can be
represented by an $n\times n$ matrix $\mathcal{S}(\mathcal{L})$ $=(\ell
_{ij}),$ called its \emph{structure matrix}\textbf{\emph{\ }} with entries
$\ell_{ij}=\ell_{i}\star\ell_{j}$ from its set $\mathcal{L}$ of elements. All
abstract properties of the groupoid $\left(  \mathcal{L},\star\right)  $ are
embodied in this matrix.

\begin{definition}
Let $\left(  \mathcal{L},\star\right)  $ be any finite groupoid of order $n$,
where $\mathcal{L}=\{\ell_{x}\mid x=1,...,n\}$ and $\star$ is a closed binary
operation on $\mathcal{L}$. The $n\times n$ matrix $\mathcal{S}(\mathcal{L})$
$=(\ell_{ij})$, where $\ell_{ij}=\ell_{i}\star\ell_{j}\in\mathcal{L}$ for all
$i,j=1,...,n$, is the \textbf{structure matrix }(or \textbf{Cayley table}) of
$\left(  \mathcal{L},\star\right)  $.\vspace{0.05in}
\end{definition}

\begin{center}%
\[
\underset{\text{(A) Structure Matric}}{\left[
\begin{array}
[c]{cccccc}%
\ell_{11} & \ell_{12} & \cdots & \ell_{1j} & \cdots & \ell_{1n}\\
\ell_{21} & \ell_{22} & \cdots & \ell_{2j} & \cdots & \ell_{2n}\\
\vdots & \vdots & \ddots & \vdots & \vdots & \vdots\\
\ell_{i1} & \ell_{i2} & \cdots & \ell_{ij} & \cdots & \ell_{in}\\
\vdots & \vdots & \vdots & \vdots & \ddots & \vdots\\
\ell_{n1} & \ell_{n2} & \cdots & \ell_{nj} & \cdots & \ell_{nn}%
\end{array}
\right]  }\qquad\underset{\text{(B) Cayley Table}}{
\begin{tabular}
[c]{|c|cccccc|}\hline
$\star$ & $\ell_{1}$ & $\ell_{2}$ & $\cdots$ & $\ell_{j}$ & $\cdots$ &
$\ell_{n}$\\\hline
$\ell_{1}$ & $\ell_{11}$ & $\ell_{12}$ & $\cdots$ & $\ell_{1j}$ & $\cdots$ &
$\ell_{1n}$\\
$\ell_{2}$ & $\ell_{21}$ & $\ell_{22}$ & $\cdots$ & $\ell_{2j}$ & $\cdots$ &
$\ell_{2n}$\\
$\vdots$ & $\vdots$ & $\vdots$ & $\ddots$ & $\vdots$ & $\vdots$ & $\vdots$\\
$\ell_{i}$ & $\ell_{i1}$ & $\ell_{i2}$ & $\cdots$ & $\ell_{ij}$ & $\cdots$ &
$\ell_{in}$\\
$\vdots$ & $\vdots$ & $\vdots$ & $\vdots$ & $\vdots$ & $\ddots$ & $\vdots$\\
$\ell_{n}$ & $\ell_{n1}$ & $\ell_{n2}$ & $\cdots$ & $\ell_{nj}$ & $\cdots$ &
$\ell_{nn}$\\\hline
\end{tabular}
}%
\]
\vspace{0.15in}
\end{center}

\begin{quote}
{\small Figure 1. General form of the} \textbf{\emph{structure matrix}%
\ }({\small or simply }\textbf{\emph{$\mathcal{S}$-matrix}}) $\mathcal{S}%
(\mathcal{L})=(\ell_{ij})${\small \ of the groupoid} $\left(  \mathcal{L}%
,\star\right)  $ {\small of order }$n${\small . The entries }$\ell_{ij}%
=\ell_{i}\star\ell_{j}${\small \ are elements of }$\mathcal{L}$. {\small This
matrix, when provided with row and column headings, is also commonly called a}
\textbf{\emph{Cayley table}}.\vspace{0.15in}
\end{quote}

Being a matrix, $\mathcal{S}(\mathcal{L})$ can be formally subjected under
certain conditions to some matrix operations like addition, transpose, etc.
The labeling of the elements of$\ \mathcal{L}$ is arbitrary; it is determined
simply by considerations of convenience and may be changed without affecting
the structure of the system. However, if $\mathcal{L}$ has an identity element
$\ell_{1}$, it is convenient to arrange its entries such that $\ell_{1k}%
=\ell_{k1}=\ell_{k}$ for all $k=1,...,n$. We shall call this arrangement the
\textbf{\emph{standard }}or \textbf{\emph{normal form}} of $\mathcal{S}%
(\mathcal{L})$ and often simplify the notation for $\mathcal{S}(\mathcal{L})$
entries by writing $i\equiv\ell_{i}$ for all $i=1,...,n$.\bigskip

\subsection{\textbf{Existence of Non-Associative Finite Invertible Loops
(NAFIL loops)\medskip}}

Now that we have defined what NAFIL loops are, our first main task is to
establish their existence. It can be shown [1] that any invertible loop of
order $n\leq4$ is a group and that if $n=5$, there exists a unique
non-associative invertible loop [3, 8]. Hence, this loop is the smallest
NAFIL. If $n\geq5,$ there are two special families of NAFIL loops whose
members cover all ODD and EVEN orders from $n=5$ to infinity [3, 6]. The
$\mathcal{S}$-matrices (Cayley tables) of the first members of each of these
families are shown in Figure 2. Thus, we have

\begin{theorem}
There exists at least one NAFIL (non-associative finite invertible loop) of
every finite order $n\geq5$.
\end{theorem}

This is the \emph{Fundamental Theorem} of the \emph{Theory of Non-Associative
Finite Invertible Loops}. From this theorem we find that the smallest NAFIL is
of order $n=5$. Moreover, it can be shown that there is exactly one (up to
isomorphism) such loop of this order whose Cayley table is that shown in
Figure 2 for $n=5.$\vspace{0.075in}

\begin{center}%
\[%
\begin{array}
[c]{ccc}%
ODD\,\left(  n=5\right)  & \qquad\qquad & EVEN\,\left(  n=6\right) \\
\left[
\begin{array}
[c]{ccccc}%
\text{1} & \text{2} & \text{3} & \text{4} & \text{5}\\
\text{2} & \text{1} & \text{5} & \text{3} & \text{4}\\
\text{3} & \text{4} & \text{1} & \text{5} & \text{2}\\
\text{4} & \text{5} & \text{2} & \text{1} & \text{3}\\
\text{5} & \text{3} & \text{4} & \text{2} & \text{1}%
\end{array}
\right]  &  & \left[
\begin{array}
[c]{cccccc}%
\text{1} & \text{2} & \text{3} & \text{4} & \text{5} & \text{6}\\
\text{2} & \text{1} & \text{4} & \text{5} & \text{6} & \text{3}\\
\text{3} & \text{4} & \text{1} & \text{6} & \text{2} & \text{5}\\
\text{4} & \text{5} & \text{6} & \text{1} & \text{3} & \text{2}\\
\text{5} & \text{6} & \text{2} & \text{3} & \text{1} & \text{4}\\
\text{6} & \text{3} & \text{5} & \text{2} & \text{4} & \text{1}%
\end{array}
\right] \\
(L_{5},\ast) &  & (L_{6},\star)
\end{array}
\]

\end{center}

\begin{quote}
{\small Figure 2. The }$\mathcal{S}${\small -matrices of the first members of
the ODD and EVEN families of NAFIL loops: (L}$_{5},\ast${\small ) and (L}%
$_{6},\star${\small ).\bigskip}
\end{quote}

\begin{remark}
\emph{From 1996 to 2001, we undertook a study [8, 15] to determine all
nonisomorphic NAFIL loops of small order n = 5, 6, and 7. This study showed
that there is exactly one NAFIL loop of order n = 5, 33 of order n = 6, and
2,333 of order n = 7. }
\end{remark}

\subsection{\textbf{Subsystem Composition of NAFIL Loops}}

Most algebraic systems contain smaller systems (called \emph{subsystems}) in
their structures [13]. Thus algebras may contain \emph{subalgebras} and groups
can have \emph{subgroups}. Similarly, certain loops (like NAFILs) may also
have \emph{subloops}.

By definition, a non-empty subset $H$ of a set $G$ is a subgroupoid (subloop,
subgroup) of a groupoid $(G,\circ)$ iff $(H,\circ)$ is a groupoid (loop,
group) [9, 13]. This means that $(H,\circ)$ satisfies all axioms satisfied by
$(G,\circ).$ In the case of invertible loops like the NAFILs, we have shown
that $m\leq\frac{n}{2},$ where $m$ is the order of $H$ and $n$ is the order of
$G.$ Not every groupoid, however, can have subsystems.

In general algebraic systems can be classified into two types: \emph{composite
}(with at least one non-trivial subsystem) and \emph{non-composite }or
\emph{plain} [14] (without any non-trivial subsystem).

\section{\textbf{Basic Concepts and Terminology}}

The NAFIL is a natural generalization of the group because it satisfies all
group axioms except A6 (associative axiom). Hence, all theorems and basic
algebraic concepts of group theory that do not depend on A6 (e.g.
\emph{homomorphisms, isomorphisms, cosets, subsystems, quotients, etc.}) apply
to NAFIL loops. Those that depend on A6 like the fundamental theorem of
Lagrange on the order of a subgroup are not satisfied in general by NAFIL
loops. Similarly, it is not always possible to uniquely define in the
traditional sense such concepts as the \emph{power} and \emph{order} of a
NAFIL element.

Because of the applications of NAFIL loops in non-associative algebras and
general loop theory, we shall also adapt several basic ideas and terms from
these fields like those involved in alternative and power-associative
algebras, Moufang, Bol, and IP loops, etc. Moreover, we shall also introduce
new concepts and ideas in addition to well established basic algebraic concepts.

In what follows, unless otherwise indicated, we shall mainly consider\emph{\ }%
NAFIL loops. And, for convenience of notation, we shall sometimes represent
$\left(  \mathcal{L},\star\right)  $ simply by $\mathcal{L}$ and the product
$x\star y$ by $xy$ (juxtaposition) if no confusion arises.

\subsection{\textbf{Products and Powers of NAFIL Elements}}

Because NAFIL loops are by definition non-associative, expressions involving
products of more than two elements are meaningful only if the grouping of the
factors is defined. Such groupings (or parenthesizings) must be clearly
indicated in all expressions. To simplify matters, it becomes necessary to
define special products of NAFIL elements.

\begin{definition}
Let $\left(  \mathcal{L},\star\right)  $ be a NAFIL and let $\ell_{1},\ell
_{2},\ell_{3},...,\ell_{m}\in\mathcal{L}$. The \textbf{left} and \textbf{right
products} of these $m$ elements are defined by the expressions:
\begin{align*}
Left  & :L(\ell_{1}\ell_{2}\ell_{3}\ldots\ell_{m})\equiv\ell_{m}(\cdots
(\ell_{3}(\ell_{2}\ell_{1}))\cdots)\\
Right  & :R(\ell_{1}\ell_{2}\ell_{3}\ldots\ell_{m})\equiv(\cdots((\ell_{1}%
\ell_{2})\ell_{3})\cdots)\ell_{m}%
\end{align*}

\end{definition}

As noted earlier, the idea of the power of a NAFIL element is not always
meaningful. However, certain special forms of repeated products of an element
by itself can be defined that are useful in many contexts. If we set
$\ell\equiv\ell_{x}$ for all $x=1,...,m$ in Definition 4, then we have:

\begin{definition}
The expression $L(\ell^{m})\equiv\underset{to\,m\,factors}{\underbrace{\text{
}\ell(\cdots(\ell(\ell\ell))\cdots)\,}},$ where $m$ is any positive integer,
is the m-th\textbf{\ left power }of $\ell\in\mathcal{L}$ and $R(\ell
^{m})\!\equiv\underset{to\,m\,factors}{\underbrace{\text{ }(\cdots((\ell
\ell)\ell)\ldots)\ell\text{ }}}$is the m-th \textbf{right power} of $\ell
\in\mathcal{L}$. If $R(\ell^{m})=L(\ell^{m}) $, we shall simply write
$\ell^{m}$ to represent either $R(\ell^{m})$ or $L(\ell^{m}). $
\end{definition}

The above definition can be stated in the equivalent form: Let $\ell
^{\lambda(m)}\equiv L(\ell^{m})$ and $\ell^{\rho(m)}\equiv R(\ell^{m}),$ where
$\lambda$ means \emph{left power} and $\rho$ means \emph{right power}.

It is often convenient to define $R(\ell^{m})$ recursively as follows: Let
$R(\ell^{1})=\ell,$ $R(\ell^{2})=\ell\ell,$ $R(\ell^{3})=R(\ell^{2})\ell,$ and
$R(\ell^{m})=R(\ell^{m-1})\ell$ whenever $m>1.$ Then we can write $R(\ell
^{m})=R(\ell^{m-1})\ell$ in the form: $\ell^{\rho(m)}=\ell^{\rho(m-1)}\ell.$
Similarly, we can also write $L(\ell^{m})=\ell L(\ell^{m-1})$ as
$\ell^{\lambda(m)}=\ell\ell^{\lambda(m-1)}.$ If $R(\ell^{m})=L(\ell^{m})$ as
in the case of abelian loops, then we simply drop the $R,$ $L,$ $\rho$ and
$\lambda$ and write: $\ell^{m}=\ell^{m-1}\ell=\ell\ell^{m-1}.$

\subsection{\textbf{Generators and Order}}

In finite group theory where the power of a group element $\ell\in\mathcal{L}$
is uniquely defined, if $n$ is the smallest positive integer such that
$\ell^{n}=1,$ where $1$ is the identity, then the elements of the set
$\mathcal{L=}\left\{  \ell^{x}\mid x=1,...,n\right\}  $ are distinct and
$\mathcal{L}$ is said to be generated by $\ell.$ If $\left(  \mathcal{L}%
,\star\right)  $ is a group$,$ any element $\ell\in\mathcal{L}$ whose powers
generate $\mathcal{L}$ is therefore called a\emph{\ generator} of the group.
Moreover, if $m\leq n$ is the smallest positive integer for which $\ell
^{m}=1,$ then $m$ is called the \emph{order} of the element $\ell
\in\mathcal{L}.$ In this case, the set $\left\{  \ell^{x}\mid
x=1,...,m\right\}  $ of powers of $\ell$ always forms a subgroup of order $m.$
This is not true in general for loops like the NAFILs where the set of powers
of an element does not necessarily form a subsystem. Nevertheless, there are
many interesting cases where this condition is satisfied.

The idea of the \emph{order} of an element and a \emph{generator} of a set can
also be extended to NAFIL loops in a modified form as follows:

\begin{definition}
Let $\left(  \mathcal{L},\star\right)  $ be a NAFIL loop of order $n$ whose
identity is $1.$ (a) The \textbf{order }of an element $\ell\in\mathcal{L}$ is
the order $m\leq n$ of the smallest subsystem generated by $\ell.$ If $m=n,$
then we say that $\ell$ is a \textbf{generator} of $\left(  \mathcal{L}%
,\star\right)  .$ (b) The \textbf{left (right) power-order} of $\ell$ is the
least positive integer $m\leq n$ such that $\ell^{\lambda(m)}=1$ $(\ell
^{\rho(m)}=1).$ If $m=n,$ then we say that $\ell$ is a \textbf{left (right)
power-generator} of $\mathcal{L}.$
\end{definition}

Definition 7(a) is the standard definition of \emph{order} which is based on
the fact that in a loop every element generates a subsystem. This definition
applies to all loops. However, unlike groups, the elements of the generated
subsystem are not restricted to the powers of an element. Because of this we
introduce in Definition 7(b) the idea of the \emph{power-order} (or simply
\emph{p-order}) of an element and that of a \emph{power-generator }(or simply
\emph{p-generator}) of a set which are analogous to the ideas of order and
generator in group theory. If $m$ is the least positive integer such that
$\ell^{\lambda(m)}=1,$ then by Definition 7(b) the order $m$ of the set
$\left\{  \ell^{\lambda(x)}\mid x=1,...,m\right\}  $ of left powers of
$\ell\ $is called the \emph{Left p-order }of $\ell$ to distinguish it from the
standard \emph{order} given in Definition 7(a): the former is simply the order
of the set of left powers of the element, whereas the latter is the order of
the smallest subsystem generated by the element. The \emph{Right p-order} of
$\ell$ is similarly defined. This distinction is important because for
non-abelian loops the Left and Right p-orders of an element are not
necessarily equal. Moreover, the set of powers of an element does not
necessarily form a subsystem. If the set of Left (Right) powers of $\ell$
forms a subsystem, then its Left (Right) p-order is equal to its standard
\emph{order}. This condition is completely satisfied by loops that are
power-associative. For such loops, there is no fundamental distinction between
order and p-order as well as between generator and p-generator.

A NAFIL loop $\left(  \mathcal{L},\star\right)  $ is called
\textbf{\emph{monogenic }}[13] if it can be generated by a single element. If
a monogenic loop can be generated by the powers of a given element $\ell
\in\mathcal{L}$, then we can distinguish two kinds of generators: \emph{left
p-generator} \emph{(LpG)} and \emph{right p-generator (RpG)} according as
$\mathcal{L}$ is generated by the left or right powers of the element $\ell$.
The distinction between \emph{left }and\emph{\ right p-generator} (as well as
between \emph{left} and \emph{right} \emph{p-order}) is useful when dealing
with loops with left- or right- handed properties.

It must be pointed out that the set $\mathcal{L}$ may not have any generator
at all other than the set itself. In many cases, it may not also have a single
element that can generate it. However, it could happen that it may be
generated by several elements that constitute a \textbf{\emph{set of
generators}}. Also, the set $\mathcal{L}$ may have more than one set of generators.

Following conventional notation, if $\ell\in\mathcal{L},$ then the set of all
elements generated by $\ell$ shall be denoted by $\left\langle \ell
\right\rangle .$ If a set is generated by two or more elements $q_{1,}%
q_{2,}...,q_{m}$, we shall denote this by $\left\langle q_{1,}q_{2,}%
...,q_{m}\right\rangle .$ However, the set $\left\langle \ell\right\rangle $
of elements generated by the element $\ell$\ is understood to mean not only
the powers of $\ell$ but may also include elements that are not powers of
$\ell$ $.$ In this case, some elements of $\left\langle \ell\right\rangle $
may be products of the powers of $\ell,$ their inverses, etc.\bigskip

\begin{remark}
\emph{The terms} \textbf{\emph{order}} \emph{of an element and}
\textbf{\emph{generator}} \emph{of a set as given in Definition 7(a) do not
indicate how they are determined in specific cases. For instance, given a loop
}$\left(  \mathcal{L},\star\right)  $\emph{\ of order }$n,$\emph{\ where
}$\mathcal{L}=\{\ell_{i}\mid i=1,2,...,n\},$\emph{\ the set of elements
generated by an element }$\ell_{i}$\emph{\ (denoted by }$\left\langle \ell
_{i}\right\rangle $\emph{) consists of all powers (left and right) of }%
$\ell_{i},$\emph{\ the products of these powers}$,$\emph{\ the products of
}$\ell_{i}$\emph{\ and its powers, etc. Computationally, this is very
difficult to determine for each element of }$\mathcal{L}.$ \emph{The
}\textbf{\emph{p-order}} \emph{concept given in Definition 7(b), on the other
hand, is easier to determine and its use is often more appropriate in most
cases than the standard }\textbf{\emph{order.}}\bigskip
\end{remark}

\subsection{\textbf{Association Properties}}

Weaker forms of A6, called \textbf{\emph{weak associative laws}}, are known to
play important roles in loop theory and many of these also apply to NAFIL
loops. These are special \emph{identities} (or \emph{identical relations})
which have the general form of the associative relation, $(ab)c=a(bc),$ but
which hold true only under certain conditions. These \emph{weak associative
laws} shall also be called \textbf{\emph{association properties}}\emph{. }

In general, a weak associative law [11] is defined as a universally quantified
equation of the form $\alpha=\beta,$ where for some variables, $V_{1}%
,V_{2},...,V_{n}$ (not necessarily distinct), $\alpha$ and $\beta$ are both
products of the form $V_{1}V_{2}...V_{n}$ (with some distribution of
parentheses). The number $n$ of variables is also called the \emph{size} of
the equation. Such a law is called \emph{non-trivial} iff $\alpha$ and $\beta$
do not have the same distribution of parentheses. For instance, the relation
$(x(yz))x=(xy)(zx),$ called a Moufang identity, is an equation of size $n=4$
that is non-trivial. Here, only three of the four variables are distinct. If
the equation $\alpha=\beta$ is such that any variable in it appears exactly
once on each side, then it is called a \emph{balanced identity}. The simplest
example of this is the associative relation: $(xy)z=x(yz).$

\subsubsection{\textbf{Inverse Properties}}

Many NAFIL loops are known to satisfy certain important identities [9]
involving weak forms of A6, called \textbf{\emph{inverse properties}}, such as
those given in

\begin{definition}
Let $\left(  \mathcal{L},\star\right)  $ be a NAFIL and let $q,\ell,\ell
^{-1}\in\mathcal{L}$. If $\ell^{-1}(\ell q)=q$ and $(q\ell)\ell^{-1}=q,$ then
$\left(  \mathcal{L},\star\right)  $ is said to have the \textbf{Left Inverse
Property (LIP)} and \textbf{Right Inverse Property (RIP)}, respectively. If
$\left(  \mathcal{L},\star\right)  $ satisfies both LIP and RIP properties,
then $\ell^{-1}(\ell q)=(q\ell)\ell^{-1}=q$ and it is said to have the
\textbf{Inverse Property (IP). }
\end{definition}

The equations that define the above inverse properties are balanced identities
of size three. Thus, the LIP identity has the form $\ell^{-1}(\ell
q)=(\ell^{-1}\ell)q=q$ since $(\ell^{-1}\ell)=1$, where $1$ is the identity
element. Similarly, the RIP identity has the form $(q\ell)\ell^{-1}=q(\ell
\ell^{-1})=q$ since $(\ell\ell^{-1})=1.$ This implies that in an IP loop,
$\ell^{-1}\ell=\ell\ell^{-1}=1$ so that every element has a unique two-sided
inverse. Therefore, all IP loops are NAFILs but the converse is not true:
there are NAFIL loops that are not IP loops.

In most loops, A3 (\emph{Inverse axiom}) is not assumed to hold. These loops
are not invertible but their elements satisfy weak forms of A3. Thus, if
$\left(  \mathcal{L},\star\right)  $ is a loop whose identity element is
$\ell_{1}\equiv1,$ then there exists elements $\ell^{-\lambda}$ and
$\ell^{-\rho},$ called the \textbf{\emph{left inverse}} and
\textbf{\emph{right inverse }} of $\ell\in\mathcal{L},$ respectively, such
that $\ell^{-\lambda}\star\ell=1$ and $\ell\star\ell^{-\rho}=1$. If $\left(
\mathcal{L},\star\right)  $ is an invertible loop, however, then
$\ell^{-\lambda}=\ell^{-\rho}\equiv\ell^{-1}$ such that $\ell^{-1}\ell
=\ell\ell^{-1}=1$ because in this case every $\ell\in\mathcal{L}$ has a unique
(or \textbf{\emph{two-sided}}) \textbf{\emph{inverse}} $\ell^{-1}%
\in\mathcal{L}.$ For loops that are not invertible, the inverse properties
given in Definition 7 will still hold provided that the weak forms of A3 are considered.

It is easy to show that in an IP NAFIL, the linear equations $ax=b$ and $ya=b
$ have the unique solutions $x=a^{-1}b$ and $y=ba^{-1},$ respectively. On the
other hand, in a LIP NAFIL only the equation $ax=b$ has the unique solution
$x=a^{-1}b$ while in a RIP NAFIL only the equation $ya=b$ has the unique
solution $y=ba^{-1}.$

\subsubsection{\textbf{Other Weak Associative Properties}}

\begin{definition}
Let $\left(  \mathcal{L},\star\right)  $ be any NAFIL and let $\ell_{i}%
,\ell_{k}\in\mathcal{L}.$ If $\left(  \mathcal{L},\star\right)  $ satisfies
the identities
\[%
\begin{tabular}
[c]{c}%
$\ell_{i}(\ell_{i}\ell_{k})=\ell_{i}^{2}\ell_{k}$
\textbf{[left\ alternative\ property\ (LAP)]}\\
$(\ell_{i}\ell_{k})\ell_{k}=\ell_{i}\ell_{k}^{2}$ \textbf{[right alternative
property (RAP)}]
\end{tabular}
\]
for all elements $\ell_{i},\ell_{k}\in\mathcal{L}$, then $\left(
\mathcal{L},\star\right)  $ is called \textbf{alternative } and is said to
have the \textbf{Alternative Property (AP)}. If $\left(  \mathcal{L}%
,\star\right)  $ satisfies the identity
\[
\ell_{i}(\ell_{k}\ell_{i})=(\ell_{i}\ell_{k})\ell_{i}%
\ \text{[\textbf{flexible\ law\ (FL)]}}%
\]
then it is called \textbf{flexible.}
\end{definition}

If $\ell_{i}=\ell_{k}=\ell$ in the LAP and RAP identities, then $L(\ell
^{3})=\ell\ell^{2}=\ell^{2}\ell=R(\ell^{3})$ and we find that $L(\ell
^{3})=R(\ell^{3}).$ In general, we find that $L(\ell^{u})=\ell\ell^{u-1}%
=\ell^{u-1}\ell=R(\ell^{u})$ so that $L(\ell^{u})=R(\ell^{u})=\ell^{u}$, where
$u$ is any positive integer. Moreover, it can be shown that an AP loop
$\left(  \mathcal{L},\star\right)  $ is power-associative. It is important to
note, however, that there are loops that satisfy only some of these
identities. Thus there are \emph{LAP loops, RAP loops,} and \emph{FL loops}.

The loops that have been studied most extensively are the
\textbf{\emph{Moufang loops}}\emph{\ }[9] because they are closest to groups
in their properties. Such loops are NAFIL loops because they are
non-associative and they satisfy all invertible loop axioms. The smallest
Moufang loop is of order $n=12$ which is non-abelian. And the smallest abelian
Moufang loop is of order $n=81.$

\begin{definition}
Let $\left(  \mathcal{L},\star\right)  $ be a NAFIL and let $\ell_{i}%
,\,\ell_{j},\,\ell_{k}\in\mathcal{L}.$ If
\[
\ell_{i}[\,\ell_{j}(\,\ell_{i}\ell_{k})]=[(\ell_{i}\,\ell_{j})\,\ell_{i}%
]\ell_{k}%
\]
then $\left(  \mathcal{L},\star\right)  $ is called a \textbf{Moufang loop}
and is said to have the \textbf{Moufang Property (MP)}.
\end{definition}

The expression
\begin{equation}
\ell_{i}[\,\ell_{j}(\,\ell_{i}\ell_{k})]=[(\ell_{i}\,\ell_{j})\,\ell_{i}%
]\ell_{k}\tag{D9.1}%
\end{equation}
is known as the \textbf{\emph{Moufang identity}} and it is equivalent to each
of the identities:
\begin{equation}
\ell_{i}[\,\ell_{j}(\,\ell_{k}\ell_{j})]=[(\ell_{i}\,\ell_{j})\,\ell_{k}%
]\ell_{j}\tag{D9.2}%
\end{equation}
\begin{equation}
(\ell_{i}\,\ell_{j})(\,\ell_{k}\ell_{i})=\ell_{i}[(\ell_{j}\,\ell_{k})\ell
_{i}]\tag{D9.3}%
\end{equation}

It can be shown that a Moufang loop satisfies the alternative property (AP)
and therefore also the power-associative property (PAP). For if we let
$\ell_{k}=1$ (identity element) in (D9.1), (D9.2) and (D9.3), then we obtain
the identities: $\ell_{i}(\,\ell_{j}\ell_{i})=(\ell_{i}\,\ell_{j})\ell_{i}$,
$\ell_{i}(\,\ell_{j}\ell_{j})=(\ell_{i}\,\ell_{j})\ell_{j}$ and $(\ell
_{i}\,\ell_{j})\ell_{i}=\ell_{i}(\ell_{j}\ell_{i}).$ These three identities
satisfy the requirements of Definition 8 for alternativity.

Another interesting structure is the Bol loop [13] which is closely related to
the Moufang loop. In a sense, it is a generalization of the Moufang loop.

\begin{definition}
Let $\left(  \mathcal{L},\star\right)  $ be a NAFIL and let $\ell_{i}%
,\,\ell_{j},\,\ell_{k}\in\mathcal{L}.$ If
\[
\left[  (\ell_{i}\ell_{j})\ell_{k}\right]  \ell_{j}=\ell_{i}\left[  (\ell
_{j}\ell_{k})\ell_{j}\right]
\]
then $\left(  \mathcal{L},\star\right)  $ is called a \textbf{right Bol loop
(RBol)}. If it satisfies the identity
\[
\left[  \ell_{i}(\ell_{j}\ell_{i})\right]  \ell_{k}=\ell_{i}\left[  \ell
_{j}(\ell_{i}\ell_{k})\right]
\]
then $\left(  \mathcal{L},\star\right)  $ is called a \textbf{left Bol loop
(LBol)}.
\end{definition}

It is known that there is a duality between the right and left Bol loops. A
right Bol loop is RIP and RAP and a left Bol loop is LIP and LAP. Given the
Cayley table of a right Bol loop, its \emph{transpose} [9] is the Cayley table
of a left Bol loop. Thus the distinction between them is not fundamental. The
smallest right (left) Bol loop is of order $n=8$; and there are exactly 6
right (left) Bol loops of this order. Because of this duality, we shall simply
call such loops as \emph{Bol loops}.

There is an interesting variety of loops, called \textbf{\emph{extra loops}},
satisfying the following equivalent identities%

\begin{align*}
(x(yz))y  & =(xy)(zy)\quad\text{and}\quad(yz)(yx)=y((zy)x)\\
((xy)z)x  & =x(y(zx))
\end{align*}
Thus, if a loop satisfies any of these identities, it will also satisfy the
others. These are weak associative laws of size $4$ in $3$ distinct variables.
If we set $x=1$ in the first two equations and $z=1$ in the in the third, we
find that $(yz)y=y(zy)$ and $(xy)x=x(yx),$ respectively. Hence, extra loops
are also flexible. Note that the first two equations are the \emph{mirrors}
[11] of each other, that is, one can be obtained from the other by writing it
backwards. A loop that satisfies any of the above identities is said to have
the \textbf{\emph{extra loop property (ELP).}}

Because the left and right powers of a NAFIL element are not always equal, the
familiar law of exponents, $\ell^{a}\star\ell^{b}=\ell^{a+b},$ is not always
satisfied. However, there are certain NAFIL loops in which this law holds.

\begin{definition}
Let $\left(  \mathcal{L},\star\right)  $ be a NAFIL and let $\ell
\in\mathcal{L}$. If $\ell^{a}\star\ell^{b}=\ell^{a+b}$, where $a$ and $b$ are
any two positive integers, then $\left(  \mathcal{L},\star\right)  $ is called
\textbf{power-associative} and is said to have the \textbf{Power Associative
Property (PAP)}.
\end{definition}

Although the PAP identity does not have the explicit form of a weak
associative law, it is a consequence of the associative axiom A6. The defining
identities are infinite in number; they are special cases of A6 involving one
variable (e.g. $x^{2}x=xx^{2},(x^{2}x)x=x^{2}x^{2}=x(xx^{2}),$...) [9]. This
indicates that in a power-associative NAFIL, A6 holds in a limited way.

\subsection{\textbf{Other Loop Properties}}

There are other loop properties defined by identities that do not have the
form of the weak associative properties. Among these are the following [9]:
\emph{Weak Inverse Property (WIP),\ Automorphic Inverse Property (AIP),
Semiautomorphic Inverse Property (SAIP), Antiautomorphic Inverse Property
(AAIP), }and the \emph{Crossed Inverse Property (CIP)}. These are defined as follows:

\begin{definition}
A NAFIL $\left(  \mathcal{L},\star\right)  $ is said to have: (a) the
\textbf{Crossed Inverse Property (CIP) }if it satisfies the identity
\[
(\ell q)\ell^{-1}=q
\]
for all $\ell,q\in\mathcal{L},$ (b) the \textbf{Weak Inverse Property (WIP)
}\emph{\ }if it satisfies the identity
\[
\ell(q\ell)^{-1}=q^{-1}%
\]
for all $\ell,q\in\mathcal{L},$ and (c) the \textbf{Automorphic Inverse
Property (AIP)}\emph{\ }if it satisfies the identity
\[
(\ell q)^{-1}=\ell^{-1}q^{-1}%
\]
for all $\ell,q\in\mathcal{L}$.
\end{definition}

There are two other known inverse properties that are related to AIP:

\begin{itemize}
\item \emph{SAIP (semiautomorphic inverse property): }There are two forms

\emph{RSAIP}: $((\ell q)\ell)^{-1}=(\ell^{-1}q^{-1})\ell^{-1}$ and
\emph{LSAIP}: $(\ell(q\ell))^{-1}=\ell^{-1}(q^{-1}\ell^{-1}).$ A loop that
satisfies both RSAIP and LSAIP is simply called SAIP.

\item \emph{AAIP (antiautomorphic inverse property): }$(\ell q)^{-1}%
=q^{-1}\ell^{-1}.\bigskip$
\end{itemize}

The above definitions of CIP, WIP, and AIP hold for certain invertible loops.
An example of a NAFIL that is CIP, WIP, and AIP is the non-abelian loop
$\left(  L_{5},\ast\right)  $ of order $5$ whose Cayley table is shown in
Figure 1. This NAFIL is also a PAP and an FL loop. Any NAFIL in which every
element is self-inverse (\textbf{\emph{unipotent}}) is trivially
power-associative. In such a loop, every element generates a subgroup of order
2 and is also called \textbf{\emph{monassociative}}. Thus, all loops belonging
to the EVEN family of NAFIL loops are monassociative.\medskip

Examples of loops that are LSAIP and RSAIP are as follows:

\begin{center}%
\[
\underset{\text{LSAIP NAFIL of order 9 }}{\underset{}{
\begin{tabular}
[c]{|c|lllllllll|}\hline
$\circ$ & 1 & 2 & 3 & 4 & 5 & 6 & 7 & 8 & 9\\\hline
1 & 1 & 2 & 3 & 4 & 5 & 6 & 7 & 8 & 9\\
2 & 2 & 3 & 4 & 1 & 6 & 7 & 8 & 9 & 5\\
3 & 3 & 4 & 1 & 2 & 7 & 8 & 9 & 5 & 6\\
4 & 4 & 1 & 2 & 3 & 8 & 9 & 5 & 6 & 7\\
5 & 5 & 6 & 7 & 9 & 1 & 2 & 3 & 4 & 8\\
6 & 6 & 7 & 8 & 5 & 9 & 1 & 2 & 3 & 4\\
7 & 7 & 8 & 9 & 6 & 4 & 5 & 1 & 3 & 3\\
8 & 8 & 9 & 5 & 7 & 3 & 4 & 6 & 1 & 2\\
9 & 9 & 5 & 6 & 8 & 2 & 3 & 4 & 7 & 1\\\hline
\end{tabular}
}}\quad\underleftrightarrow{transpose}\quad\underset{\text{RSAIP NAFIL of
order 9}}{\underset{}{
\begin{tabular}
[c]{|c|ccccccccc|}\hline
$\circ^{\prime}$ & 1 & 2 & 3 & 4 & 5 & 6 & 7 & 8 & 9\\\hline
1 & 1 & 2 & 3 & 4 & 5 & 6 & 7 & 8 & 9\\
2 & 2 & 3 & 4 & 1 & 6 & 7 & 8 & 9 & 5\\
3 & 3 & 4 & 1 & 2 & 7 & 8 & 9 & 5 & 6\\
4 & 4 & 1 & 2 & 3 & 9 & 5 & 6 & 7 & 8\\
5 & 5 & 6 & 7 & 8 & 1 & 9 & 4 & 3 & 2\\
6 & 6 & 7 & 8 & 9 & 2 & 1 & 5 & 4 & 3\\
7 & 7 & 8 & 9 & 5 & 3 & 2 & 1 & 6 & 4\\
8 & 8 & 9 & 5 & 6 & 4 & 3 & 2 & 1 & 7\\
9 & 9 & 5 & 6 & 7 & 8 & 4 & 3 & 2 & 1\\\hline
\end{tabular}
}}%
\]
\bigskip
\end{center}

These loops are the transpose of each other. Both of them are PAP and have the
same subsystems: \{1,2,3,4\}, \{1,3\}, \{1,5\}, \{1,6\}, \{1,7\}, \{1,8\}, \{1,9\}.

Another type of loop is called \textbf{\emph{totally symmetric (TS)}} (also
called a \emph{Steiner loop}) which satisfies the identities%

\[
\ell_{x}\ell_{y}=\ell_{y}\ell_{x}\text{\quad and\quad}\ell_{x}(\ell_{x}%
\ell_{y})=\ell_{y}%
\]
If such a loop is a NAFIL, we shall call it a \textbf{\emph{TS NAFIL}}. Here
again, the equation $\ell_{x}(\ell_{x}\ell_{y})=\ell_{y}$ has the form
$\ell_{x}(\ell_{x}\ell_{y})=(\ell_{x}\ell_{x})\ell_{y}=\ell_{y},$ where
$\ell_{x}\ell_{x}=1.$ It is clear from this that a TS NAFIL is an abelian IP
loop such that every element is self-inverse.

\subsubsection{\textbf{Summary of Loop Properties}}

Other weak associative laws have also been found useful in the study of loops
and quasigroups. However, many of these do not apply to NAFIL loops nor do
they contribute much to the understanding of these structures.

A concise summary of various loop properties known to be satisfied by NAFILs
is given in the table below. These properties are defined by universally
quantified equations called \emph{identities} or \emph{identical relations}.

\bigskip

\begin{center}
$\underset{\text{Table 3. List of known special loop properties used to test
NAFIL loops of orders }n\,=\,5,\,6,\,7.}{\underset{}{
\begin{tabular}
[c]{|c|c|c|}\hline
\textbf{Special Loop Property} & \textbf{Acronym} & \textbf{Defining
Equation}\\\hline\hline
\multicolumn{1}{|l|}{{\small Left Inverse Property}} &
\multicolumn{1}{|l|}{{\small LIP}} & $\ell^{-1}(\ell q)=(\ell^{-1}\ell
)q=q$\\\hline
\multicolumn{1}{|l|}{{\small Right Inverse Property}} &
\multicolumn{1}{|l|}{{\small RIP}} & $(q\ell)\ell^{-1}=q(\ell\ell^{-1}%
)=q$\\\hline
\multicolumn{1}{|l|}{{\small Inverse Property}} &
\multicolumn{1}{|l|}{{\small IP}} & {\small LIP and RIP}\\\hline
\multicolumn{1}{|l|}{{\small Left Alternative Property}} &
\multicolumn{1}{|l|}{{\small LAP}} & $\ell_{i}(\ell_{i}\ell_{k})=(\ell_{i}%
\ell_{i})\ell_{k}$\\\hline
\multicolumn{1}{|l|}{{\small Right Alternative Property}} &
\multicolumn{1}{|l|}{{\small RAP}} & $(\ell_{i}\ell_{k})\ell_{k}=\ell_{i}%
(\ell_{k}\ell_{k})$\\\hline
\multicolumn{1}{|l|}{{\small Alternative Property}} &
\multicolumn{1}{|l|}{{\small AP}} & {\small LAP and RAP}\\\hline
\multicolumn{1}{|l|}{{\small Flexible Law}} & \multicolumn{1}{|l|}{{\small FL}%
} & $\ell_{i}(\ell_{k}\ell_{i})=(\ell_{i}\ell_{k})\ell_{i}$\\\hline
\multicolumn{1}{|l|}{{\small Moufang Property}} &
\multicolumn{1}{|l|}{{\small MP}} & $\ell_{i}[\ell_{j}(\ell_{i}\ell
_{k})]=[(\ell_{i}\ell_{j})\ell_{i}]\ell_{k}$\\\hline
\multicolumn{1}{|l|}{{\small Left Bol}} & \multicolumn{1}{|l|}{{\small LBol}}
& $\left[  \ell_{i}(\ell_{j}\ell_{i})\right]  \ell_{k}=\ell_{i}\left[
\ell_{j}(\ell_{i}\ell_{k})\right]  $\\\hline
\multicolumn{1}{|l|}{{\small Right Bol}} & \multicolumn{1}{|l|}{{\small RBol}}
& $\left[  (\ell_{i}\ell_{j})\ell_{k}\right]  \ell_{j}=\ell_{i}\left[
(\ell_{j}\ell_{k})\ell_{j}\right]  $\\\hline
\multicolumn{1}{|l|}{{\small Extra Loop Property}} &
\multicolumn{1}{|l|}{{\small ELP}} & $(\ell_{i}(\ell_{j}\ell_{k}))\ell
_{j}=(\ell_{i}\ell_{j})(\ell_{k}\ell_{j})$\\\hline
\multicolumn{1}{|l|}{{\small C Loop Property}} &
\multicolumn{1}{|l|}{{\small CP}} & $\ell_{i}[\ell_{j}(\ell_{j}\ell
_{k})]=[(\ell_{i}\ell_{j})\ell_{j}]\ell_{k}$\\\hline
\multicolumn{1}{|l|}{{\small RIF Loop Property}} &
\multicolumn{1}{|l|}{{\small RIFP}} & \multicolumn{1}{|l|}{$(\ell_{i}\ell
_{j})[\ell_{k}(\ell_{i}\ell_{j})]=\{[\ell_{i}(\ell_{j}\ell_{k})]\ell_{i}%
\}\ell_{j}$}\\\hline
\multicolumn{1}{|l|}{{\small A} {\small sub m} {\small Loop Property}} &
\multicolumn{1}{|l|}{{\small A\_m} {\small P}} & \multicolumn{1}{|l|}{$\ell
_{i}[(\ell_{j}\ell_{i})(\ell_{k}\ell_{i})]=[(\ell_{i}\ell_{j})(\ell_{i}%
\ell_{k})]\ell_{i}$}\\\hline
&  & \\\hline
\multicolumn{1}{|l|}{{\small Power Associative Property}} & {\small PAP} &
$\ell^{a}\star\ell^{b}=\ell^{a+b}$\\\hline
\multicolumn{1}{|l|}{{\small Totally Symmetric}} & {\small TS} & $\ell_{i}%
\ell_{j}=\ell_{j}\ell_{i}\text{\ and\ }\ell_{i}(\ell_{i}\ell_{j})=\ell_{j}%
$\\\hline
\multicolumn{1}{|l|}{{\small Weak Inverse Property}} &
\multicolumn{1}{|l|}{{\small WIP}} & $\ell(q\ell)^{-1}=q^{-1}$\\\hline
\multicolumn{1}{|l|}{{\small Automorphic Inverse Property}} &
\multicolumn{1}{|l|}{{\small AIP}} & $(\ell q)^{-1}=\ell^{-1}q^{-1}$\\\hline
\multicolumn{1}{|l|}{{\small Anti-Automorphic Inverse Property}} &
\multicolumn{1}{|l|}{{\small AAIP}} & $(\ell q)^{-1}=q^{-1}\ell^{-1}$\\\hline
\multicolumn{1}{|l|}{{\small Left Semi-Automorphic Inverse Prop.}} &
\multicolumn{1}{|l|}{{\small LSAIP}} & $((\ell q)\ell)^{-1}=(\ell^{-1}%
q^{-1})\ell^{-1}$\\\hline
\multicolumn{1}{|l|}{{\small Right Semi-Automorphic Inverse Prop.}} &
\multicolumn{1}{|l|}{{\small RSAIP}} & $(\ell(q\ell))^{-1}=\ell^{-1}%
(q^{-1}\ell^{-1})$\\\hline
\multicolumn{1}{|l|}{{\small Semi-Automorphic Inverse Property}} &
\multicolumn{1}{|l|}{{\small SAIP}} & {\small LSAIP and RSAIP}\\\hline
\multicolumn{1}{|l|}{{\small Crossed Inverse Property}} &
\multicolumn{1}{|l|}{{\small CIP}} & $(\ell q)\ell^{-1}=q$\\\hline
\multicolumn{1}{|l|}{{\small Left Cojugacy Closed}} &
\multicolumn{1}{|l|}{{\small LCC}} & $\ell_{i}(\ell_{j}\ell_{k})=[\ell
_{i}(\ell_{j}\ell_{i}^{-1})](\ell_{i}\ell_{k})$\\\hline
\multicolumn{1}{|l|}{{\small Right Cojugacy Closed}} &
\multicolumn{1}{|l|}{{\small RCC}} & $(\ell_{i}\ell_{j)}\ell_{k}=(\ell_{i}%
\ell_{k})[(\ell_{k}^{-1}\ell_{j})\ell_{k}]$\\\hline
\multicolumn{1}{|l|}{{\small Conjugacy Closed Loop Property}} &
\multicolumn{1}{|l|}{{\small CCP}} & {\small LCC and RCC}\\\hline
\end{tabular}
}}$
\end{center}

\section{\textbf{Construction and Analysis of NAFIL Loops}}

So far, we have defined the NAFIL and presented some of its basic properties.
Our next important question is: How do we construct such a loop? Finally, how
do we analyze the constructed loop to determine its properties?

To answer these questions we must first specify the kind of loop we propose to
construct. Usually, a groupoid (loop, group, quasigroup) is constructed by
means of a set of independent generators and a set of relations. However, this
procedure is often difficult to carry out. We mentioned earlier that a finite
groupoid is completely defined by its structure matrix or Cayley table. We
know that the structure matrix of a quasigroup is a Latin square. Since a loop
is a quasigroup with an identity element, then its Cayley table is a Latin
Square in normal (or reduced) form. Thus, the construction of a loop (like the
NAFIL or group) is equivalent to the construction of a Latin square. This
method of construction is highly developed and it can be done manually or by
means of computer programs [5].

Because there are many kinds of NAFILs, there is no general method of
constructing them. A survey of the literature has shown that there are
numerous ad hoc methods of construction [9]. We have, however, developed a
simple and efficient method of construction and analysis of finite loops
called the \emph{Structure Matrix Method}. This method can be carried out
manually as well as by a computer program called \emph{FINITAS}.[7]

\bigskip


\bigskip

\begin{thebibliography}{99}                                                                                               %
\bibitem {ref1}A.A. Albert, \emph{Quasigroups. I}, Trans. Amer. Math. Soc. 54
(1943), 507-519. [See also: .\emph{Quasigroups. II}, Trans. Amer. Math. Soc.
55 (1944), 401-419.]

\bibitem {ref2}E. T. Bell, \emph{The Development of Mathematics}, McGraw-Hill
Book Company, New York (1945).

\bibitem {ref3}R. E. Cawagas, \emph{On the Existence of Pseudogroups},
Matimyas Matematika, \textbf{5} (1981), 8-13.

\bibitem {ref4}\_\_\_\_\textit{\ }\emph{Construction of All Cayley algebras of
Order 2}$^{r}$\emph{\ by the $\mathcal{Z}$SM Process}, Abstracts of Short
Communications, International Congress of Mathematicians, Zurich (1994), 7.
[See also: R. E. Cawagas, Transactions, National Academy of Science and
Technology, XV (1993), 133-142.]

\bibitem {ref5}\_\_\_\_ \emph{Computer-Based Construction and Analysis of
Finite Algebras}, Proceedings of the International Conference on Computational
Mathematics, Chulalongkorn University, Thailand (1997), 107-115.

\bibitem {ref6}\_\_\_\_ \emph{Foundations of the Theory of Finite
Pseudogroups}, PUP SciTech R\&D\ Center, Manila (1996). [See also: R. E.
Cawagas, \emph{Foundations of the Theory of Finite Pseudogroups},
Transactions, National Academy of Science and Technology, XVI (1994), 175]

\bibitem {ref7}\_\_\_\_ \emph{FINITAS -- A Software for the Construction and
Analysis of Finite Algebraic Structures}, PUP Journal of Research and
Exposition (1997), Vol. 1, No. 1 , pp. 1-10. [See also: R. E. Cawagas,
\emph{AXIOMS - Software for the Construction and Analysis of Finite
Quasigroups, Semigroups and Related Structures,} Proceedings of the Second
Asian Mathematical Conference 1995, World Scientific, Singapore-New
Jersey-London-Hong Kong (1998), 401 405.]

\bibitem {ref8}\_\_\_\_ \emph{Generation of NAFIL loops of small order},
Quasigroups and Related Systems 7(2000), pp. 1-5.

\bibitem {ref9}O. Chein, et al (Editors) \emph{Quasigroups and Loops: Theory
and Applications}, Sigma Series in Pure Mathematics, Helderman Verlag Berlin (1990).

\bibitem {ref10}D. Hobby and R. McKenzie, \emph{Structure of Finite Algebras},
American Mathematical Society (1988).

\bibitem {}K. Kunen, Quasigroups, Loops, and Associative Laws, J. Algebra 185
(1996) 194-204.

\bibitem {}S. Okubo, \emph{Introduction to Octonion and Other Non-Associative
Algebras in Physics}, Cambridge University Press, Cambridge (1995).

\bibitem {}H. O. Pflugfelder, \emph{Quasigroups and Loops: Introduction},
Sigma Series in Pure Mathematics, Helderman Verlag Berlin (1990).

\bibitem {}J. D. H. Smith, \emph{Mal'cev Varieties}, Lecture Notes in
Mathematics 554, Springer-Verlag, Berlin$\cdot$Heidelberg$\cdot$New York
(1976), pp. 96-112..

\bibitem {ref11}J. Zhang and H. Zhang, \emph{Generating models by SEM}, Proc.
of International Conference on Automated Deduction (CADE-96), pp. 308-312.
[See also: H. Zhang, et al, \emph{PSATO: a distributed propositional prover
and its application to quasigroup problems}, Journal of Symbolic Computation
(1996) 21, 543-560.]
\end{thebibliography}
\end{document}